\documentclass[12pt,twoside,reqno]{amsart}
\usepackage{amssymb,verbatim}
\usepackage{color}
\usepackage{hyperref}
\theoremstyle{definition}

\usepackage[english]{babel}
\usepackage{amsmath}
\usepackage{amsfonts}
\usepackage{blindtext}
\usepackage{tabularx}
\usepackage{tabulary}
\usepackage{multirow}
\usepackage{graphicx}
\usepackage[all]{xy}
\usepackage{tikz}
\usepackage{amsthm}
\usepackage{enumerate}
\usepackage{titlesec}
\usepackage{lipsum}
\usepackage{bbm}

\titleformat{\section}
  {\normalfont\scshape\filcenter}{\thesection}{1em}{}

\numberwithin{equation}{section}

\titleformat{\subsection}
{\normalfont\normalsize\bfseries}{\thesubsection.}{1em}{}

\theoremstyle{plain}
	\newtheorem{Thm}{Theorem}[section]
	\newtheorem{Lem}[Thm]{Lemma}
	\newtheorem{Cor}[Thm]{Corollary}
	
	\newtheorem{Prop}[Thm]{Proposition}
	\newtheorem{Rem}[Thm]{Remark}

\theoremstyle{definition}
	
	\newtheorem{Ex}[Thm]{Example}

\theoremstyle{remark}

\begin{document}
\title[Hochschild Cohomology]{The Hochschild Cohomology \\ of square-free monomial complete intersections}

\author{Nghia T.H. Tran}
\address{School of Mathematics, Statistics and Applied Mathematics, National University of Ireland, Galway, Ireland}
\email{\{n.tranthihieu1, emil.skoldberg\}@nuigalway.ie}
\urladdr{}

\author{Emil Sk\" oldberg}
\address{}
\email{}
\urladdr{}

\thanks{Version: \today}

\maketitle
\begin{abstract}
We determine the Hochschild cohomology algebras of the square-free monomial complete intersections. In particular we provide a formula for the cup product which gives the cohomology module an algebra structure and then we provide a description of this structure in terms of generators and relations. In addition, we compute the Hilbert series of the Hochschild cohomology of these algebras.
\end{abstract}
\section*{Introduction}
The Hochschild (co)homology of an algebra contains substantial information about the algebra. Many works have studied different aspects of the structure of Hochschild (co)homology (see \cite{Ger63, Gro91, Hol96, CS97}). However, the multiplication in the Hochschild (co)homology is complicated to calculate and not known in general. In \cite{Hol00}, T. Holm determined the ring structure of the Hochschild cohomology for the $k$-algebra $k[X]/\langle f \rangle$ where $f$ is a monic polynomial in the polynomial ring $k[X]$ of a single variable $X$ and $k$ is a commutative ring with unity. In particular, this description was given in terms of generators and relations, which is fruitful to utilize in some sense.\\

Our work was initially inspired by attempts to understand the structure of Hochschild cohomology of algebras with more variables. We consider the square-free monomial intersections, which is in the form of algebras $A=k\left[x_1,x_2,\ldots,x_n\right]/\langle m_1,m_2,\ldots,m_r\rangle$ where $m_i$'s are square-free such that $\mathrm{supp}(m_i)\cap \mathrm{supp}(m_j)=\emptyset$. After relabeling, we have $$A\cong \dfrac{k\left[x_1,\ldots,x_{u_1}\right]}{x_1\cdots x_{u_1}}\otimes_k \cdots \otimes_k \dfrac{k\left[x_1,\ldots,x_{u_s}\right]}{x_1\cdots x_{u_s}}\otimes_k k\left[x_1,\ldots,x_m\right].$$
As the tensor product is preserved under the action of taking the Hochschild cohomology, i.e., $\mathrm{
HH}^*(A\otimes_k B)\cong \mathrm{HH}^*(A)\otimes_k\mathrm{HH}^*(B)$, we only need to study the algebras of the form  $k\left[x_1,x_2,\ldots,x_n\right]/\langle x_1x_2\cdots x_n\rangle$ where $k$ is a field and $n \in \mathbb{N}$.\\

The note is divided into six sections. In Section \ref{s.1}, we review some basic notions of Hochschild cohomology. Also, here we interpret the alternative resolution given by Guccione et al. in \cite{GGRV92} to construct the Hochschild cohomology module in our case. In Section \ref{s.2}, we study the $k$-space structure of the  Hochschild cohomology of the given algebra, denoted by $\mathrm{HH}^*(A)$. In particular, we describe this module via the cohomology of simpler complexes in Theorem \ref{theorem1}. In Sections \ref{s.3} and \ref{s.4}, starting from a cocycle we construct an explicit chain map between the shifted resolution and the resolution itself. From this map, we can describe the Yoneda product which gives the $k$-space an algebra structure. In Section \ref{s.5}, we use the previous results to describe the ring structure of $\mathrm{HH}^*(A)$ in terms of generators and relations. This is the main result of the paper, which is presented in Theorem \ref{main.theorem}.
In Section \ref{s.6}, we finish the article by computing the Hilbert series of $\mathrm{HH}^*(A)$, which is given in Theorem \ref{hilbert.theorem}.\\

\textbf{Acknowledgements.} We are grateful to some of our colleagues for their careful reading and useful comments on the paper. Our warmest thanks go to Isaac Burke for discussions on the Hochschild (co)homology topic. In order to check the computations in examples, we used Macaulay2 \cite{GS}.
\section{A construction of Hochschild cohomology}\label{s.1}
We first introduce some notation. Let $A^{op}$ be the opposite algebra of $A$. In our case, $A^{op}$ is exactly $A$ since $A$ is commutative. Throughout this paper, tensor products will be taken over $k$, i.e., $\otimes=\otimes_k$ unless otherwise indicated. By $A^e=A\otimes A^{op}$ we denote the enveloping algebra of $A$. We use the same notation for elements in $k[x_1,x_2,\ldots,x_n]$ and their cosets in the quotient ring $A$ with the convention $x_1x_2 \cdots x_n=0$, if there are no ambiguities.\\

In order to compute the Hochschild cohomology, we use the resolution in \cite{GGRV92} which is more conducive to explicit computation than the historical resolution of Hochschild, which was given in \cite{Hoc45}, (see also \cite{Wit17}). We give details of the resolution in case of the square-free monomial $f=x_1x_2\cdots x_n$ and get the following resolution:
\begin{equation} \label{p1complex1}
\cdots \xrightarrow{\text{ }\text{ }\text{ }\text{$d_3$}\text{ }\text{ }\text{ }} F_2 \xrightarrow{\text{ }\text{ }\text{ }\text{$d_2$}\text{ }\text{ }\text{ }} F_1 \xrightarrow{\text{\text{ }\text{ }\text{ }$d_1$}\text{ }\text{ }\text{ }} F_0\xrightarrow{\text{\text{ }\text{ }\text{ }$\epsilon$}\text{ }\text{ }\text{ }} A\longrightarrow 0
\end{equation}
where $F_m$ is a finitely generated free $A^e$-module with basis elements $e_{i_1\cdots i_s}\cdot t^{(q)}$ ($s$, $q$ $\geq 0$ and $s+2q=m$), whereby $e_{i_1\cdots i_s}$ or $e_I$ with $I=\{i_1,\ldots, i_s\}$ we mean $e_{i_1}\wedge \cdots \wedge e_{i_s}$ ($1 \leq i_1 <\cdots <i_s \leq n$). Then \eqref{p1complex1} is an exact sequence of free $A^e$-modules $F_m$ with
$$\epsilon: A^e \to A,\qquad a\otimes b \mapsto ab$$
and where the differentials $d_m$ are defined inductively as follows:
\begin{align*}\label{f1}
d_s(e_{i_1\cdots i_s})&=
\sum\limits_{j=1}^s(-1)^{j-1}(1\otimes x_{i_j}-x_{i_j}\otimes 1)e_{i_1\cdots \widehat{i}_j \cdots i_s};\\
d_2(t)&=\sum\limits_{j=1}^n x_1\cdots x_{j-1}\otimes x_{j+1}	\cdots x_n\cdot e_j;\\
d_{s+2q}(e_{i_1\cdots i_s} t^{(q)})&= d_s(e_{i_1\cdots i_s}) t^{(q)}+d_2(t)\cdot e_{i_1\cdots i_s} \cdot t^{(q-1)}, \text{ if } q \geq 1.
\end{align*}
For abbreviation, we sometimes write $d$ instead of $d_m$.\\

Applying the contravariant functor $\mathrm{Hom}_{A^e}(-,A)$ to the truncation of the above resolution, we obtain a complex of $A^e$-modules and $A^e$-homomorphisms
\begin{equation*}
0 \longrightarrow \mathrm{Hom}_{A^e}(F_0,A)  \xrightarrow {\text{ }\text{ }\text{ }d^*_1\text{ }\text{ }\text{ }} \mathrm{Hom}_{A^e}(F_1,A) \xrightarrow {\text{ }\text{ }\text{ }d^*_2\text{ }\text{ }\text{ }}  \mathrm{Hom}_{A^e}(F_2,A) \longrightarrow \cdots.
\end{equation*}
From the last complex, by passing to cohomology one gets the Hochschild cohomology $\mathrm{HH}^*(A)$ of $A$. In more detail,
\begin{equation*}
\mathrm{HH}^*(A)=\mathop \oplus \limits_{i\geq 0}\mathrm{HH}^i(A),
\end{equation*}
that is, $\mathrm{HH}^i(A)=\mathrm{H}^i(\mathrm{Hom}_{A^e}(F_\bullet,A))=\mathrm{Ker}(d^*_{i+1})/\mathrm{Im}(d^*_{i})$ for all $i \geq 0$, where $d^*_{0}$ is taken to be the zero map.\\

So far we have constructed the Hochschild cohomology of algebra $A$ to be an $\mathbb{N}$-graded $A^e$-module.

\section{Hochschild cohomology $\mathrm{HH}^*(A)$ as a $k$-space}\label{s.2}
In this section, we consider the Hochschild cohomology as a graded $k$-space and give a description of the structure of this module via simpler complexes.\\

For any $m \in \mathbb{N}$, let $\overline{F}_m$ be the $k$-space spanned by the same basis elements as $F_m$. By the definition of $F_m$, the number of basis elements is finite. There is an isomorphism between the following $k$-spaces
\begin{equation*}
\mathrm{Hom}_{A^e}(F_m,A) \cong  \mathrm{Hom}_k(\overline{F}_m,A)
\end{equation*}
for all $m \in \mathbb{N}$.
Thus, we get a new complex of $k$-spaces and $k$-homomorphisms
\begin{equation} \label{complex2}
0 \longrightarrow \text{Hom}_{k}(\overline{F}_0,A)  \xrightarrow{\text{ }\text{ }\text{ }\partial_1\text{ }\text{ }\text{ }} \text{Hom}_{k}(\overline{F}_1,A) \xrightarrow{\text{ }\text{ }\text{ }\partial_2\text{ }\text{ }\text{ }} \text{Hom}_{k}(\overline{F}_2,A) \longrightarrow\cdots,
\end{equation}
where the maps $\partial_m$ (for $m \in \mathbb{N}$) will be stated in the subsequent lemma. For abbreviation, we often use $\partial$ instead of $\partial_m$ when we do not need to specify the index $m$.\\

Now let us introduce some notation which will appear in the sequel:
\begin{itemize}
\item $[n]:=\{1,2,\ldots,n\}$;
\item $\mathrm{sgn}(i, I):=(-1)^{\left|\left\{j \in I\mid j<i\right\}\right|}$ where $|S|$ is the cardinality of the set $S$;
\item $\mathbf{x}^\alpha:=x_1^{\alpha_1}x_2^{\alpha_2}\cdots x_n^{\alpha_n}$ where $\alpha$ is the lattice point $(\alpha_1,\alpha_2,\ldots, \alpha_n)$ in $\mathbb{N}^n$;
\item $\mathrm{supp}(\mathbf{x}^\alpha):=\{i\mid \alpha_i>0\}$ where  $\alpha=(\alpha_1,\alpha_2,\ldots, \alpha_n)$ as above.
\end{itemize}
Let $e_It^{(q)}$ be a basis element in $\overline{F}_m$ and $\mathbf{x}^\alpha$ be a basis element in $A$. We denote  by $(e_It^{(q)},\mathbf{x}^\alpha)$ the $k$-linear map in $\text{Hom}_{k}(\overline{F}_m,A)$ which sends $e_It^{(q)}$ to $\mathbf{x}^\alpha$ and other basis elements to 0, i.e.,
\begin{equation*}
(e_It^{(q)},\mathbf{x}^\alpha)(e_Jt^{(p)}) =
  \begin{cases}
    \mathbf{x}^\alpha,       & \quad \text{if } J=I \text{ and }p=q,\\
    0,  & \quad  \text{otherwise}.\\
  \end{cases}
\end{equation*}
Let us call these basis elements the \textit{standard} elements; since $\overline{F}_m$ is finite dimensional, they form a basis of $\mathrm{Hom}_k(\overline{F}_m,A)$. We also use the same notation, $(e_It^{(q)},\mathbf{x}^\alpha)$, for the residue class in $\mathrm{HH}^*(A)$.

\begin{Lem} \label{lemma1}
Let $(e_It^{(q)},\mathbf{x}^\alpha)$ be a standard element in $\mathrm{Hom}_{k}(\overline{F}_m,A)$. We have
\begin{equation*}
\partial_{m+1}(e_It^{(q)},\mathbf{x}^\alpha)=\sum\limits_{i \in I\smallsetminus \mathrm{supp}(\mathbf{x}^\alpha)}{\mathrm{sgn}(i,I)(e_{I\smallsetminus \{i\}}t^{(q+1)},\mathbf{x}^\alpha \cdot x_1\cdots x_{i-1}\cdot x_{i+1}\cdots x_n)}.
\end{equation*}
\proof
Let us consider the following diagram
\begin{equation*}
\xymatrix{\mathrm{Hom}_{A^e}(F_m,A) \ar[d]_{\cong} \ar[r]^{d^*_{m+1}} & \mathrm{Hom}_{A^e}(F_{m+1},A)\ar[d]_{\cong}\\
\mathrm{Hom}_{k}(\overline{F}_m,A) \ar[u] \ar[r]^{\partial_{m+1}}&\mathrm{Hom}_{k}(\overline{F}_{m+1},A). \ar[u]}
\end{equation*}
By combining the isomorphisms, we can deduce $\partial_{m+1}$ from $d^*_{m+1}$ straightforwardly.\\
As $f=(e_It^{(q)},\mathbf{x}^\alpha)$ is a standard element of $\mathrm{Hom}_{k}(\overline{F}_m,A)$, $f$ is identified with a function in $\mathrm{Hom}_{A^e}(F_m,A)$. A direct calculation shows that it is $(e_It^{(q)},\mathbf{x}^\alpha)$, which is also denoted by $f$ by abuse of notation. The canonical homomorphism
\begin{equation*}
\begin{array}{cccc}
d^*_{m+1}: & \mathrm{Hom}_{A^e}(F_m,A) & \longrightarrow & \mathrm{Hom}_{A^e}(F_{m+1}, A)\\
       & f                & \longmapsto &d^*_{m+1}(f):=f\circ d
\end{array}
\end{equation*}
sends $f$ to $d^*_{m+1}(f)$, an $A^e$-homomorphism from $F_{m+1}$ to $A$.\\

Recall that $A$ is a left and right $A^e$-module by the scalar multiplication: $(a\otimes b)\cdot c=acb$ and $c\cdot(a\otimes b) =bca$ respectively.
Let  $e_{i_1\cdots i_s}t^{(p)}$(shortly, $e_Jt^{(p)}$) be a basis element in $F_{m+1}$. By \eqref{f1}, one has
\begin{equation*}
d^*_{m+1}(f)(e_Jt^{(p)})=f\left(d_{m+1}(e_Jt^{(p)})\right)= A_1+A_2
\end{equation*}
where
\begin{align*}
A_1&=f\left( \sum\limits_{j=1}^{s}(-1)^{j+1}d_{1}(e_{i_j}) e_{J\smallsetminus \{i_j\}} t^{(p)}\right)=\sum\limits_{j=1}^{s}(-1)^{j+1}(1\otimes x_{i_j}-x_{i_j}\otimes 1)\cdot f( e_{J\smallsetminus \{i_j\}} t^{(p)})\\
&=\sum\limits_{j=1}^{s}(-1)^{j+1}\left(f( e_{J\smallsetminus \{i_j\}} t^{(p)})\cdot x_{i_j}-x_{i_j}\cdot f( e_{J\smallsetminus \{i_j\}} t^{(p)})\right)
=0; \text{ and}
\end{align*}
\begin{align*}
A_2&=f\left((-1)^{l}e_J d_2(t) t^{(p-1)}\right)=
f\left( \sum\limits_{i=1}^n{x_1\cdots x_{i-1}\otimes x_{i+1}\cdots x_n\cdot e_i\wedge e_J t^{(p-1)}}\right)\\
&=\sum\limits_{i\in [n]\smallsetminus J}{\mathrm{sgn}(i,J) x_1\cdots x_{i-1}\cdot f(e_{J\cup \{i\}}\cdot t^{(p-1)})}\cdot x_{i+1}\cdots x_{n}\\
&=\sum\limits_{i\in [n]\smallsetminus J}{\mathrm{sgn}(i,J) x_1\cdots x_{i-1}\cdot x_{i+1}\cdots x_n\cdot f(e_{J\cup \{i\}}\cdot t^{(p-1)})}.
\end{align*}
Since
\begin{equation*}
f(e_{J\cup \{i\}}\cdot t^{(p-1)})=
\begin{cases}
\mathbf{x}^\alpha & \text{ if } J\cup \{i\} = I \text{ and } p-1=q\\
0 & \text{otherwise}
\end{cases},
\end{equation*}
we have
\begin{equation*}
A_2= \begin{cases}
\mathrm{sgn}(i,J) x_1\cdots x_{i-1}\cdot x_{i+1}\cdots x_n\cdot\mathbf{x}^\alpha & \text{ if } J = I\smallsetminus \{i\}  \text{ and } p=q+1\\
0 & \text{otherwise}
\end{cases}
\end{equation*}
where $i$ is an element in $I$. Note that $\mathrm{sgn}(i,J)=\mathrm{sgn}(i,I)$ when $J = I\smallsetminus \{i\}$ and if $i \in \mathrm{supp}(\mathbf{x}^\alpha)$, then $x_1\cdots x_{i-1}\cdot x_{i+1}\cdots x_n\cdot\mathbf{x}^\alpha=0$ in $A$. Hence, the formula follows.
\qed
\end{Lem}
From this lemma, we are going to derive some consequences about properties of standard elements in the kernel and the image of $\partial$.
\begin{Cor} \label{cor1.2}
Let $(e_It^{(q)},\mathbf{x}^\alpha)$ be a standard element. We have $\partial (e_It^{(q)},\mathbf{x}^\alpha)=0$ if and only if $I$ is a subset of $\mathrm{supp}(\mathbf{x}^\alpha)$.
\end{Cor}
\begin{Cor} \label{lemma2}
The non-zero standard element $(e_It^{(q)},\mathbf{x}^\alpha)$ occurs as a component of some element in $\mathrm{Im}(\partial)$ if and only if the two following conditions hold:
\begin{enumerate} [(i)]
	\item $q>0$; and
	\item There exists some index $i$ in $ [n]\smallsetminus I$ which satisfies $\mathrm{supp}(\mathbf{x}^\alpha)=[n]\smallsetminus \{i\}$.
\end{enumerate}
\proof The first part $``\Rightarrow"$ is straightforward by observing the formula of $\partial$ in Lemma \ref{lemma1}. Conversely, if $(e_It^{(q)},\mathbf{x}^\alpha)$ satisfies the two conditions (i) and (ii), it will be a component in the image of the element $\left(e_{I\cup\{i\}}t^{(q-1)},\dfrac{\mathbf{x}^\alpha}{x_1\cdots \widehat{x_i}\cdots x_n}\right)$.
\qed
\end{Cor}
\begin{Cor} \label{cor1}
If $(e_It^{(q)},\mathbf{x}^\alpha)$ and $(e_Jt^{(p)},\mathbf{x}^\beta)$ are standard elements such that their images under $\partial$ have some non-zero component in common, then they are identical.
\proof This is an immediate consequence of Corollary \ref{lemma2}. \qed
\end{Cor}

The $k$-space $\mathrm{Hom}_k(\overline{F}_m,A)$ is generated by standard elements of degree $m$, i.e.,
\begin{equation*}
\mathrm{Hom}_k(\overline{F}_m,A)=\mathop \oplus \limits_{|I|+2q=m}k(e_It^{(q)},\mathbf{x}^\alpha)
\end{equation*}
where $k(e_It^{(q)},\mathbf{x}^\alpha)$ is the $k$-module generated by the standard element $(e_It^{(q)},\mathbf{x}^\alpha)$.\\
Let $\Gamma$ be the set of standard elements that are not in any component of $\mathrm{Im}(\partial)$. Let us take an element $\gamma=(e_It^{(q)},\mathbf{x}^\alpha)$ in $\Gamma$. The image of $\gamma$ under $\partial$ is given by Lemma \ref{lemma1}:
\begin{equation*}
 \partial(\gamma)=\sum\limits_{i \in I\smallsetminus \mathrm{supp}(\mathbf{x}^\alpha)}{\mathrm{sgn}(i,I)(e_{I\smallsetminus \{i\}}t^{(q+1)},\mathbf{x}^\alpha \cdot x_1\cdots x_{i-1}\cdot x_{i+1}\cdots x_n)}.
\end{equation*}
From here, we construct a complex $M_\gamma$:
\begin{equation*}
\cdots 0 \longrightarrow k(e_It^{(q)},\mathbf{x}^\alpha) \longrightarrow \mathop  \oplus \limits_{i \in I\smallsetminus \mathrm{supp}(\mathbf{x}^\alpha)}k(e_{I\smallsetminus \{i\}}t^{(q+1)},\mathbf{x}^\alpha \cdot x_1\cdots x_{i-1}\cdot x_{i+1}\cdots x_n) \longrightarrow 0 \cdots,
\end{equation*}
where the $k$-maps are taken to be the $k$-maps $\partial$ in \eqref{complex2} restricted to the corresponding subspaces. We obtain that each  such complex is a subcomplex of \eqref{complex2}, moreover it is a direct summand of \eqref{complex2}. By Corollary \ref{cor1}, the subcomplexes indexed by elements in $\Gamma$ have zero intersection. The following theorem shows that the complex \eqref{complex2} can be written as a direct sum of subcomplexes indexed by the elements in $\Gamma$.
\begin{Thm} \label{theorem1}
We have the following direct sum:
\begin{equation*}
\mathrm{Hom}_k(\overline{F}_\bullet,A)=\mathop \oplus \limits_{\gamma \in \Gamma}M_\gamma.
\end{equation*}
\proof It is obvious that we have the inclusion $\mathop \oplus \limits_{\gamma \in \Gamma}M_\gamma \subseteq \mathrm{Hom}_k(\overline{F}_\bullet,A)$ since $M_\gamma$ is a subcomplex of $\mathrm{Hom}_k(\overline{F}_\bullet,A)$ for all $\gamma \in \Gamma$. For the inverse inclusion, let us consider an arbitrary non-zero basis element $E=(e_It^{(q)},\mathbf{x}^\alpha)$ in $\mathrm{Hom}_k(\overline{F}_m,A)$. There are two conceivable cases:\\
\textbf{Case 1.} If $E$ is a component in the image $\mathrm{Im}(\partial)$, then there exists a unique element $\gamma$ in $\mathrm{Hom}_k(\overline{F}_{m-1},A)$, stated in Lemma \ref{lemma2}, not in the kernel of $\partial$ such that $\partial(\gamma)$ contains $E$ as a component. It is obvious that $\gamma \in \Gamma$ and the subcomplex $M_{\gamma}$ includes $E$.\\
\textbf{Case 2.} If $E$ is not any component in $\mathrm{Im}(\partial)$, then $E$ belongs to the subcomplex indexed by $E$ itself, $M_E$.\\
Hence, the complex \eqref{complex2} can be split into a direct sum of simpler subcomplexes as desired.
\qed
\end{Thm}
\begin{Cor} \label{cor8}
For each $i$ in $\mathbb{N}$, we have the following isomorphism:
\begin{equation*}
H^i(\mathrm{Hom}_k(\overline{F}_\bullet,A)) \cong \mathop\oplus_{\gamma \in \Gamma}H^i(M_{\gamma}).
\end{equation*}
\end{Cor}
\begin{Ex} We will see here a slice of splitting complex for the case $n=2$, which is also used to illustrate the results throughout this article. Let $A=k[x,y]/\langle xy \rangle$ where $k$ is a field. Then some of the subcomplexes in Corollary \ref{cor8} for $A$ are shown as below:
\begin{align*}
&\ldots \\
0 \longrightarrow k(e_1,y^3) \longrightarrow k(t,y^4) \longrightarrow  0\\
0 \longrightarrow k(e_1 &t,y) \longrightarrow k(t^2,y^2) \longrightarrow 0\\
0 \longrightarrow k(t,x^5) \longrightarrow 0\\
k(e_1&t,x) \\
0 \longrightarrow k(e_1e_2,1) \longrightarrow \text{ }\text{ }\text{ }\text{ }\text{ } &\oplus \text{ }\text{ }\text{ }\text{ }\text{ } \longrightarrow 0\\
k(e_2&t,y) \\
&\ldots
\end{align*}
\end{Ex}
From Corollary \ref{cor1.2}, \ref{lemma2} and Theorem \ref{theorem1}, we mention here a consequence about the kernel and image of the map $\partial$.
\begin{Rem}
The kernel of $\partial$ is spanned by the standard elements $(e_It^{(q)},\mathbf{x}^\alpha)$ where $I$ is a subset of $\mathrm{supp}(\mathbf{x}^\alpha)$. The image of $\partial$ is spanned by $\partial(e_It^{(q)},\mathbf{x}^\alpha)$ where  $I$ is not a subset of $\mathrm{supp}(\mathbf{x}^\alpha)$.
\end{Rem}
We have obtained a description of $\mathrm{HH}^*(A)$ as a $k$-module via simpler complexes. Next we will equip $\mathrm{HH}^*(A)$ with a multiplication which gives this module the structure of a $k$-algebra.
\section{An explicit chain map}\label{s.3}
For now, let $F_\bullet$ be a free resolution of the $A^e$-module $A$ and $f:F_i \to A$ be an $A^e$-homomorphism such that $f\circ d_{i+1}=0$.  By the comparison theorem, there is a chain map $\tilde f$ consisting of homomorphisms $\tilde f_m$, $m \in \mathbb{N}$ that makes the following diagram commute, moreover such a chain map is unique up to chain homotopy.
\begin{equation} \label{diagram1}
 \xymatrix{
 \cdots \ar[r]& {F_{i+2}}\ar[d]^{\tilde f_2}  \ar[r]^{d_{i+2}} &F_{i+1}\ar[d]^{\tilde f_1}\ar[r]^{d_{i+1}} &F_i\ar[d]^{\tilde f_0} \ar[rd]^f\\
\cdots \ar[r]& {F_2}  \ar[r]^{d_2} &F_1\ar[r]^{d_1} &F_0\ar[r]^{\epsilon}&A}
 \end{equation}
In theory, we know how to construct such a chain map. However, it is not generally easy in practice. The goal of this section is to provide an explicit chain map in case of our resolution that makes the above diagram commute. First we need some following auxiliary results on computations. The lemma below can be seen as a generalization of Corollary \ref{cor1.2}.
\begin{Lem} \label{lem1.chain_map}
Let $f:F_i \to A$ be a cocycle and $e_{i_1\cdots i_m}t^{(q)}$ be a basis element in $F_i$. We then have that $x_{i_j}$ is a divisor of $f(e_{i_1\cdots i_m}t^{(q)})$ for all $j \in [m]$ if $f(e_{i_1\cdots i_m}t^{(q)})\neq 0$.
\proof
Without loss of generality, we assume that $j=1$. Let us consider the element $e_{i_2\cdots i_m}t^{(q+1)} \in F_{i+1}$. Applying the differential map, one has that:
\begin{align*}
d_{i+1}(e_{i_2\cdots i_m}t^{(q+1)})=&\sum\limits_{j=2}^{m}(-1)^{j}(1\otimes x_{i_{j+1}}-x_{i_{j+1}}\otimes 1)e_{i_2\cdots \widehat{i_j}\cdots i_m}t^{(q+1)}\\
&+\sum\limits_{j=1}^n(x_1\cdots x_{j-1}\otimes x_{j+1}\cdots x_n)e_j\wedge e_{i_2\cdots i_m} t^{(q)}.
\end{align*}
Hence, $$f\circ d_{i+1}(e_{i_2\cdots i_m}t^{(q+1)})=\sum\limits_{j=1}^n(x_1\cdots x_{j-1}\cdot x_{j+1}\cdots x_n)f(e_j \wedge e_{i_2\cdots i_m} t^{(q)}).$$
Since $f\circ d_{i+1}(e_{i_2\cdots i_m}t^{(q+1)})=0$, we have that: $$\sum\limits_{j=1}^n(x_1\cdots x_{j-1}\cdot x_{j+1}\cdots x_n)f( e_j \wedge e_{i_2\cdots i_m}t^{(q)})=0.$$
Multiplying both sides by $x_1\cdots x_{i_1-1}\cdot x_{i_1+1}\cdots x_n$,  we obtain that:
$$x_1^2\cdots x^2_{i_1-1}\cdot x^2_{i_1+1}\cdots x_n^2f(e_{i_1}\wedge e_{i_2\cdots i_m}t^{(q)})=0.$$
So we must have $x_{i_1}| f(e_{i_1\cdots i_m}t^{(q)})$.
\qed
\end{Lem}
\begin{Lem} \label{lem2.chainmap}
For $z \in [n]$, set $U_z:=\sum\limits_{j=z+1}^n(x_1 \cdots x_{j-1}\otimes x_{j+1}\cdots x_n)e_j$ with the convention that $U_n=0$. We then have $d(U_z)=x_1 \cdots x_z\otimes x_{z+1}\cdots x_n.$
\proof Applying the differential map, one gets
\begin{align*}
d(U_z)=&\sum\limits_{j=z+1}^n(x_1\cdots x_{j-1}\otimes x_{j+1}\cdots x_n)(1\otimes x_j-x_j\otimes 1)\\
=&\sum\limits_{j=z+1}^n(x_1\cdots x_{j-1}\otimes x_{j}\cdots x_n-x_1\cdots x_{j}\otimes x_{j+1}\cdots x_n)
=x_1\cdots x_z\otimes x_{z+1}\cdots x_n.
\end{align*}
Thus the assertion follows.
\qed
\end{Lem}
Now we are in the position to obtain the formula of the chain map.
\begin{Prop} \label{th.chainmap}
Let $f:F_i \to A$ be a cocycle in $\mathrm{Hom}_{A^e}(F_i,A)$. For a given $j\in \mathbb{N}$, we define an $A^e$-homomorphism $\tilde f_j:F_{i+j} \to F_j$ as follows. Let $x=e_{i_1\cdots i_m}t^{(q)}$ be a basis element in $F_{i+j}$ and define:
\begin{align*}
\tilde f_j(x)=\sum\limits_{\mathcal{M}}(-1) ^{ms+j_1+\cdots +j_r-r}U_{l_s}\cdots U_{l_1}\cdot e_{i_{j_r}\cdots i_{j_1}}t^{(u)}\cdot \dfrac{f(e_{i_1\cdots \widehat{i}_{j_1}\cdots \widehat{i}_{j_r}\cdots i_m}\wedge e_{l_1\cdots l_s}t^{(q-u-s)})}{x_{l_1}\cdots x_{l_s}}\otimes 1,
\end{align*}
where the sum is indexed on $\mathcal{M}$ which consists of triples $(u,J,L)$ where $J=\{j_1,\ldots,j_r\}$, $L=\{l_1,\ldots,l_s\}$ satisfying the following conditions:
$$r+s+2u=j; 1\leq j_1 <\cdots <j_r \leq m; \text{ and } 1\leq l_1 <\cdots <l_s \leq n.$$
The chain map $\tilde f$ given as above makes the diagram \eqref{diagram1} commute.
\end{Prop}
The proof of this theorem is given by the combination of the following remarks and lemmas.

\begin{Rem} By Lemma \ref{lem1.chain_map}, we have $\dfrac{f(e_{i_1\cdots \widehat{i}_{j_1}\cdots \widehat{i}_{j_r}\cdots i_m}\wedge e_{l_1\cdots l_s}t^{(q-u-s)})}{x_{l_1}\cdots x_{l_s}}$ is an element in $A$. Also by this lemma, all elements in form of a `fraction' like above are in $A$ throughout this section.
\end{Rem}

In the our first lemma, we can see how the first step, on $\tilde f_0$, works with any homomorphism $f$, not necessarily satisfying $f\circ d_{i+1}=0$.
\begin{Lem}
For any homomorphism $f$ in $\mathrm{Hom}_{A^e}(F_i,A)$, we always have the commutative diagram below:
		\begin{equation*}
		\xymatrix{F_i \ar[d]_{\tilde f_0}\ar[rd]^f\\
F_0 \ar[r]^{\epsilon}&A}
		\end{equation*}
\proof
Indeed, for any basis element $x=e_{i_1\cdots i_m}t^{(q)}$ in $F_i$, we have $r+s+2u=0$. Then, $r=s=u=0$ is the only option and one gets that
$$\tilde f_0(x)=f(x)\otimes 1$$
So we obtain that $\epsilon  \tilde f_0=f.$
\qed
\end{Lem}
We now turn to the rest of the diagram in the following lemma.
\begin{Lem}
The homomorphisms defined in Proposition \ref{th.chainmap} make the following diagram commute:
\begin{equation*}
\xymatrix{F_{i+j} \ar[d]_{\tilde f_{j}}\ar[r]^{d_{i+j}}&F_{i+j-1} \ar[d]^{\tilde f_{j-1}}\\
F_{j} \ar[r]^{d_{j}}&F_{j-1}}
\end{equation*}
\proof
Let $f=(e_Kt^{(v)},\mathbf{x}^\alpha)$ be a standard element and $x=e_{i_1\cdots i_m}t^{(q)}$ be a basis element in $F_{i+j}$. To simplify the proof, let us introduce some notation: \begin{itemize}
\item Let $I:=\{i_1,\ldots,i_m\}$; $J:=\{i_{j_1},\ldots,i_{j_r}\}$; and $L:=\{l_1,\ldots,l_s\}$. We see that
\begin{equation*}
f(e_{i_1\cdots \widehat{i}_{j_1}\cdots \widehat{i}_{j_r}\cdots i_m}\wedge e_{l_1\cdots l_s}t^{(q-u-s)})=f(e_{(I\smallsetminus J)\cup L}t^{(q-u-s)})
\end{equation*}
up to sign, which is non-zero ($=\mathbf{x}^\alpha$) if and only if
$(I\smallsetminus J)\cup L=K$ and $q-u-s =v$.

\item $s_N=|K\smallsetminus N|$;
\item Let $M$ and $N$ be two sets of natural numbers such that $M\cap N=\emptyset$. We denote $\mathrm{sgn}(M,N)$ the power of $-1$ such that $$e_M\wedge e_N=\mathrm{sgn}(M,N)e_{M\cup N}.$$
\end{itemize}
It is simple to show that
\begin{equation*}
\mathrm{sgn}(M,N)=\prod\limits_{i \in M}{\mathrm{sgn}(i,N)}.
\end{equation*}
Let $N:=I\smallsetminus J=K\smallsetminus L$. Then we have $N \subseteq I\cap K$ and the formula becomes:
\begin{equation}\label{multi}
\tilde{f}_j(x)=\sum_{\substack{N \subseteq I\cap K\\ |N|\geq |K|-q+v}}{(\dagger)U_{K\smallsetminus N}\cdot e_{I\smallsetminus N}\cdot t^{(q-v-|K|+|N|)}\cdot \dfrac{\mathbf{x}^\alpha}{\mathbf{x}_{K\smallsetminus N}}\otimes 1},
\end{equation}
where the sign of the summand corresponding to $N$ in \ref{multi} is
\begin{equation*}
(\dagger)=(-1)^{m\cdot s_N}\cdot \mathrm{sgn}(I\smallsetminus N,N)\cdot \mathrm{sgn}(N,K\smallsetminus N);
\end{equation*}
$\mathbf{x}_M:=\prod\limits_{i \in M}x_i$; and whereby $U_M$ (for any $M=\{i_1,\ldots , i_r\}$, $i_1<\cdots <i_r$), we mean $U_{i_r}\wedge \cdots  \wedge U_{i_1}$.\\
We will show that $ d_{j}\circ \tilde f_{j} (x)= \tilde f_{j-1} \circ d_{i+j} (x)$. All the below computations are to be considered as equations up to sign. The sign of the formula will be considered later.
\begin{align*}
d_{j}\circ \tilde f_{j} (x)= & \sum_{\substack{
N \subseteq I\cap K\\
|N|\geq |K|-q+v\\
i \in K\smallsetminus N}}
{d(U_i)\cdot U_{(K\smallsetminus N)\smallsetminus \{i\}}\cdot e_{I\smallsetminus N}\cdot t^{(q-v-|K|+|N|)}\cdot \dfrac{\mathbf{x}^\alpha}{\mathbf{x}_{K\smallsetminus N}}\otimes 1}\\
&+\sum_{\substack{
N \subseteq I\cap K\\
|N|\geq |K|-q+v\\
i \in I\smallsetminus N}}
{d(e_i)\cdot U_{K\smallsetminus N}\cdot e_{(I\smallsetminus N)\smallsetminus \{i\}}\cdot t^{(q-v-|K|+|N|)}\cdot \dfrac{\mathbf{x}^\alpha}{\mathbf{x}_{K\smallsetminus N}}\otimes 1}\\
&+\sum_{\substack{
N \subseteq I\cap K\\
|N|> |K|-q+v\\
i \not\in I\smallsetminus N}}
{X_i\cdot U_{K\smallsetminus N}\cdot e_{(I\smallsetminus N)\cup \{i\}}\cdot t^{(q-v-|K|+|N|-1)}\cdot \dfrac{\mathbf{x}^\alpha}{\mathbf{x}_{K\smallsetminus N}}\otimes 1},
\end{align*}
where $X_i:= x_1\cdots x_{i-1}\otimes x_{i+1}\cdots x_n$.
Let us denote these three sums as $(1_L)$, $(2_L)$ and $(3_L)$ respectively. Computing the right hand side, we have
$$\tilde f_{j-1} \circ d_{i+j} (x)=(1_R)+(2_R)+(3_R),$$
where
\begin{align*}
(1_R):=&
\sum_{\substack{
N \subseteq I\cap K\\
|N|\geq |K|-q+v\\
i \in I\smallsetminus K}}
{d(e_i)\cdot U_{K\smallsetminus N}\cdot e_{(I\smallsetminus \{i\})\smallsetminus N}\cdot t^{(q-v-|K|+|N|)}\cdot \dfrac{\mathbf{x}^\alpha}{\mathbf{x}_{K\smallsetminus N}}\otimes 1},\\
(2_R):=&\sum_{\substack{
N \subseteq I\cap K\\
|N|\geq |K|-q+v\\
i \in (I\cap K)\smallsetminus N}}
{d(e_i)\cdot U_{K\smallsetminus N}\cdot e_{(I\smallsetminus \{i\})\smallsetminus N}\cdot t^{(q-v-|K|+|N|)}\cdot \dfrac{\mathbf{x}^\alpha}{\mathbf{x}_{K\smallsetminus N}}\otimes 1},\\
(3_R):=&[q>0]\sum_{\substack{
N \subseteq (I\cup \{i\})\cap K\\
|N|\geq |K|-q+v+1\\
i \notin I}}
{X_i\cdot U_{K\smallsetminus N}\cdot e_{(I\cup \{i\})\smallsetminus N}\cdot t^{(q-v-|K|+|N|-1)}\cdot \dfrac{\mathbf{x}^\alpha}{\mathbf{x}_{K\smallsetminus N}}\otimes 1}.
\end{align*}
Since $K\smallsetminus N=(K\smallsetminus I)\cup ((K \cap I)\smallsetminus N)$, we can write
\begin{equation*}
(1_L)=\sum\limits_{i\in K\smallsetminus N}{\ldots}=\sum\limits_{i\in K\smallsetminus I}{\ldots}+\sum\limits_{i\in (I \cap K)\smallsetminus N}{\ldots}=:(1_L)^A+(1_L)^B.
\end{equation*}
\\
Similarly, we have
\begin{equation*}
(3_L)=\sum\limits_{i\notin I\smallsetminus N}{\ldots}=\sum\limits_{i\notin I\cup K}{\ldots}+\sum\limits_{i\in N}{\ldots}+\sum\limits_{i\in K\smallsetminus I}{\ldots}=:(3_L)^A+(3_L)^B+(3_L)^C
\end{equation*}
and
\begin{equation*}
(3_R)=\sum\limits_{i\notin I}{\ldots}=\sum\limits_{i\notin I\cup K}{\ldots}+\sum\limits_{i\in (K\smallsetminus I)\cap N}{\ldots}+\sum\limits_{i\in (K\smallsetminus I)\smallsetminus N}{\ldots}=:(3_R)^A+(3_R)^B+(3_R)^C.
\end{equation*}
In order to get 
\begin{equation*}
(1_L)+(2_L)+(3_L)=(1_R)+(2_R)+(3_R),
\end{equation*}
we will show that
\begin{enumerate}[(i)]
\item $(2_L)=(1_R)+(2_R)$;
\item $(3_L)^A=(3_R)^A$;
\item $(1_L)^A=(3_R)^B$;
\item $(3_L)^B=(3_R)^C$;
\item $(1_L)^B+(3_L)^C=0$.
\end{enumerate}
As mentioned before, we will show that these sums are identical up to sign first. For (i), since $i \in I\smallsetminus N=(I\smallsetminus K)\cup ((I \cap K)\smallsetminus N)$, we have $(I\smallsetminus N)\smallsetminus \{i\} =(I\smallsetminus \{i\})\smallsetminus N$. Thus,  $(2_L)=(1_R)+(2_R)$.\\
Since $i \notin I\cup K$, $I\cap K=(I\cup \{i\})\cap K$ and $(I\smallsetminus N)\cup \{i\}=(I\cup \{i\})\smallsetminus N$. So we have (ii). \\
For (iii), Let $N':=N\cup \{i\}$ and note that $$d(U_i)\cdot \dfrac{\mathbf{x}^\alpha}{\mathbf{x}_{K\smallsetminus N}}\otimes 1=X_i\cdot \dfrac{\mathbf{x}^\alpha}{\mathbf{x}_{(K\smallsetminus N)\smallsetminus \{i\}}}\otimes 1$$ for all $i \in K\smallsetminus N$. Then we have
\begin{align*}
(1_L)^A=\sum_{\substack{
N' \subseteq (I\cup \{i\})\cap K\\
|N'|\geq |K|-q+v+1\\
i \in (K\smallsetminus I)\cap N'}}
{X_i\cdot U_{K\smallsetminus N'}\cdot e_{(I\cup \{i\})\smallsetminus N'}\cdot t^{(q-v-|K|+|N'|-1)}\cdot \dfrac{\mathbf{x}^\alpha}{\mathbf{x}_{K\smallsetminus N'}}\otimes 1},
\end{align*}
which is exactly the sum $(3_R)^B$.\\
For (iv), $N\subseteq (I\cap K)\cup \{i\}$ becomes $N\subseteq (I\cap K)$ and $(I\cup\{i\})\smallsetminus N= (I \smallsetminus N)\cup \{i\}$ for all $i \in (K\smallsetminus I)\smallsetminus N$. Then (iv) follows.
For (v), let $N'=N\smallsetminus \{i\}$. Then one has
\begin{equation*}
 (3_L)^C=\sum_{\substack{
N' \subseteq I\cap K\\
|N'|\geq |K|-q+v\\
i \in (I\cap K)\smallsetminus N'}}
{d(U_i)\cdot U_{(K\smallsetminus N')\smallsetminus \{i\}}\cdot e_{I\smallsetminus N'}\cdot t^{(q-v-|K|+|N'|)}\cdot \dfrac{\mathbf{x}^\alpha}{\mathbf{x}_{K\smallsetminus N'}}\otimes 1},
\end{equation*}
which is equal to $(1_L)^B$.\\
To complete the proof, we show that the signs of the formulas coincide. We have the sign of the summand corresponding to pair $\{i,N\}$ of the left hand side and the right hand side as follows:
\begin{align*}
\mathrm{sign}(1_L)_{\{i,N\}}&=(-1)^{m\cdot s_N}\cdot \mathrm{sgn}(I\smallsetminus N,N)\cdot \mathrm{sgn}(N,K\smallsetminus N)\cdot \mathrm{sgn}(i,K\smallsetminus N)\cdot (-1)^{s_N-1};\\
\mathrm{sign}(2_L)_{\{i,N\}}&= \mathrm{sign}(3_L)_{\{i,N\}}\\
&=(-1)^{m\cdot s_N}\cdot \mathrm{sgn}(I\smallsetminus N,N)\cdot \mathrm{sgn}(N,K\smallsetminus N)\cdot \mathrm{sgn}(i,I\smallsetminus N)\cdot (-1)^{s_N};\\
\mathrm{sign}(1_R)_{\{i,N\}}&=\mathrm{sign}(2_R)_{\{i,N\}}\\
&=(-1)^{(m-1)\cdot s_N}\cdot \mathrm{sgn}((I\smallsetminus \{i\})\smallsetminus N,N)\cdot \mathrm{sgn}(N,K\smallsetminus N)\cdot \mathrm{sgn}(i,I);\\
\mathrm{sign}(3_R)_{\{i,N\}}&=(-1)^{(m+1)\cdot s_N}\cdot \mathrm{sgn}((I\cup \{i\})\smallsetminus N,N)\cdot \mathrm{sgn}(N,K\smallsetminus N)\cdot \mathrm{sgn}(i,I).
\end{align*}
Now we are in the position to prove that the signs in equations (i) to (v) coincide.\\
 For (i), we need to show that $\mathrm{sign}(2_L)_{\{i,N\}}=\mathrm{sign}(1_R)_{\{i,N\}}$, i.e., we must have that
\begin{equation*}
\mathrm{sgn}(I\smallsetminus N,N)\cdot \mathrm{sgn}(i,I\smallsetminus N)=\mathrm{sgn}((I\smallsetminus \{i\})\smallsetminus N,N)\cdot \mathrm{sgn}(i,I).
\end{equation*}
Indeed, since $N \subseteq I$ and $i \in I\smallsetminus N$, we have $$\mathrm{sgn}(I\smallsetminus N,N)=\mathrm{sgn}((I\smallsetminus \{i\})\smallsetminus N,N)\cdot \mathrm{sgn}(i,N)$$ and  $$\mathrm{sgn}(i,N)\cdot \mathrm{sgn}(i,I\smallsetminus N)= \mathrm{sgn}(i,I).$$ Then the result follows.\\
For (ii), since $i\notin I\cup K$ and $N \subseteq I\cap K$, we have
$$\mathrm{sgn}((I\cup \{i\})\smallsetminus N,N)=\mathrm{sgn}(I\smallsetminus N,N)\cdot \mathrm{sgn}(i,N)$$ and $$ \mathrm{sgn}(i,N) \cdot\mathrm{sgn}(i,I)= \mathrm{sgn}(i,I\smallsetminus N).$$
Thus, $$\mathrm{sgn}(I\smallsetminus N,N)\cdot\mathrm{sgn}(i,I\smallsetminus N)=\mathrm{sgn}((I\cup \{i\})\smallsetminus N,N)\cdot  \mathrm{sgn}(i,I),$$
which implies that $\mathrm{sign}(3_L)^A_{\{i,N\}}=\mathrm{sign}(3_R)^A_{\{i,N\}}$.\\
For (iii), we need to show that $\mathrm{sign}(1_L)^A_{\{i,N'\}}=\mathrm{sign}(3_R)^B_{\{i,N\}}$. First we need to deduce $\mathrm{sign}(1_L)^A_{\{i,N'\}}$ from $\mathrm{sign}(1_L)^A_{\{i,N\}}$ where $N'=N\cup\{i\}$, $i \in K\smallsetminus I$, and $N \subseteq I\cap K$. We have the following identities:
\[
s_N=s_{N'}+1,
\]
and
\begin{align*}
\mathrm{sgn}((I\cup \{i\})\smallsetminus N',N')&=\prod_{\substack{
j\in (I\cup \{i\})\smallsetminus N'\\
j<i}}
{\mathrm{sgn}(j,N')}\cdot \prod_{\substack{
j\in (I\cup \{i\})\smallsetminus N'\\
j>i}}
{\mathrm{sgn}(j,N')}\\
&=\prod_{\substack{
j\in I\smallsetminus N\\
j<i}}
{\mathrm{sgn}(j,N)}\cdot \prod_{\substack{
j\in I\smallsetminus N\\
j>i}}
{\mathrm{sgn}(j,N)}\cdot (-1)^{|\{j\in I\smallsetminus N \mid j>i\}|}\\
&=\mathrm{sgn}(I\smallsetminus N,N)\cdot (-1)^{m+1-|N'|}\cdot \mathrm{sgn}(i,I\smallsetminus N').
\end{align*}
This implies that
\begin{equation*}
\mathrm{sgn}(I\smallsetminus N,N)=(-1)^{m+1-|N'|}\cdot\mathrm{sgn}((I\cup \{i\})\smallsetminus N',N')\cdot  \mathrm{sgn}(i,I\smallsetminus N').
\end{equation*}
By a similar argument, we have
\begin{equation*}
\mathrm{sgn}(N,K\smallsetminus N)=(-1)^{|N'|-1}\mathrm{sgn}(N',K\smallsetminus N')\cdot \mathrm{sgn}(i,K\smallsetminus N')\cdot \mathrm{sgn}(i,N').
\end{equation*}
Together with $\mathrm{sgn}(i,K\smallsetminus N)=\mathrm{sgn}(i,K\smallsetminus N')$, we have
\begin{align*}
\mathrm{sign}(1_L)^A_{\{i,N'\}}=(-1)^{(m+1)\cdot s_{N'}}\cdot \mathrm{sgn}((I\cup \{i\})&\smallsetminus N',N')\cdot \\
&\cdot \mathrm{sgn}(i,I\smallsetminus N')\cdot \mathrm{sgn}(N',K\smallsetminus N')\cdot \mathrm{sgn}(i,N')\\
=(-1)^{(m+1)\cdot s_{N'}}\cdot \mathrm{sgn}((I\cup \{i\})&\smallsetminus N',N')\cdot \mathrm{sgn}(N',K\smallsetminus N')\cdot \mathrm{sgn}(i,I),
\end{align*}
which is exactly $\mathrm{sign}(3_R)^B_{\{i,N\}}$ when $N'=N$.\\
Since $$\mathrm{sgn}((I\cup \{i\})\smallsetminus N,N)=\mathrm{sgn}(I\smallsetminus N,N)\cdot \mathrm{sgn}(i,N)$$ and $$\mathrm{sgn}(i,N)\cdot \mathrm{sgn}(i,I)=\mathrm{sgn}(i,I\smallsetminus N),$$ we have $$\mathrm{sgn}(I\smallsetminus N,N)\cdot \mathrm{sgn}(i,I\smallsetminus N)=\mathrm{sgn}((I\cup \{i\})\smallsetminus N,N)\cdot \mathrm{sgn}(i,I).$$ Hence, the sign in (iv) follows.\\
For the last item, (v), we will show that $\mathrm{sign}(1_B)^B_{\{i,N\}}=-\mathrm{sign}(3_L)^C_{\{i,N'\}}$. Since $i\in N \subseteq (I \cap K)$ and $N'=N\smallsetminus \{i\}$, using the same argument as item (iii), one has the following observations:
\begin{itemize}
\item $s_N=s_{N'}-1$;
\item $\mathrm{sgn}(I\smallsetminus N,N)=(-1)^{m-|N'|-1}\cdot \mathrm{sgn}(I\smallsetminus N',N')\cdot \mathrm{sgn}(i,N')\cdot\mathrm{sgn}(i,I\smallsetminus N')$;
\item $\mathrm{sgn}(N,K\smallsetminus N)=(-1)^{|N'|}\cdot \mathrm{sgn}(N',K\smallsetminus N')\cdot \mathrm{sgn}(i,K\smallsetminus N')\cdot \mathrm{sgn}(i,N')$; and
\item $\mathrm{sgn}(i,K\smallsetminus N)=\mathrm{sgn}(i,K\smallsetminus N')$.
\end{itemize}
Thus, we get $$\mathrm{sign}(3_L)^C_{\{i,N'\}}=(-1)^{m\cdot s_{N'}}\cdot \mathrm{sgn}(I\smallsetminus N',N')\cdot \mathrm{sgn}(N',K\smallsetminus N')\cdot \mathrm{sgn}(i,K\smallsetminus N')\cdot (-1)^{s_{N'}},$$
which is $-\mathrm{sign}(1_B)^B_{\{i,N\}}$ when $N'=N$.
Hence, we have the equation as desired.
\qed
\end{Lem}
\section{The cup product}\label{s.4}
In this section, we interpret the product, which is defined at the chain level on the resolution as a composition of chain maps. This product is the so-called Yoneda product and equivalent to the cup product which turns $k$-module $\mathrm{HH}^*(A)$ into a $k$-algebra (see \cite{Wit17}, Chapter 1 for more detail).\\
Let $f \in \mathrm{Hom}_{A^e}(F_i,A)$ and $g=\mathrm{Hom}_{A^e}(F_j,A)$ be cocycles. For any projective resolution $F_\bullet$ of the $A^e$-module $A$, we extend $f$ to a chain map $\tilde f:F_\bullet \to F_\bullet$ as shown in Section 3. The Yoneda product $f\smile g$ of  the two cocycles $f$ and $g$ is the composition $g\circ \tilde{f}_j$ in $\mathrm{Hom}_{A^e}(F_{i+j},A)$:
\begin{equation*}
f\smile g:=g\circ \tilde f_j.
\end{equation*}
The composition $f\smile g$ is again a cocycle. Since $\tilde f$ is unique up to homotopy, the cup product induces a well-defined product on cohomology. We will now give the product of the two standard cocycles in the following theorem, using the formula for $\tilde{f}$.
\begin{Prop} \label{th.multiplication}
Let $f=\left(e_It^{(p)},\mathbf{x}^\alpha\right)$ and $g=\left(e_Jt^{(q)},\mathbf{x}^\beta\right)$ be cocycles in $\mathrm{HH}^*(A)$. Then
\begin{equation*}
f\smile  g=
\begin{cases}
\mathrm{sgn}(I,J)\cdot \left(e_{I\cup J}t^{(p+q)},\mathbf{x}^{\alpha+\beta }\right), & \text{if } I\cap J=\emptyset, \\
0, & \text{otherwise}.
\end{cases}
\end{equation*}
Moreover, this multiplication is commutative up to sign.
\proof
Let $x=e_Kt^{(h)}$ ($K=\{i_1,\ldots,i_m\}$) be a basis element of degree $i+j$. By Proposition \ref{th.chainmap}, we have
\begin{align*}
g\circ \tilde f_j(x)=\sum\limits_{\mathcal{M}}(-1) ^{ms+j_1+\cdots +j_r-r}g(U_{l_s}\cdots U_{l_1}e_{i_{j_r}\cdots i_{j_1}}t^{(u)})\cdot \dfrac{f(e_{i_1\cdots \widehat{i}_{j_1}\cdots \widehat{i}_{j_r}\cdots i_m}\wedge e_{l_1\cdots l_s}t^{(h-u-s)})}{x_{l_1}\cdots x_{l_s}}.
\end{align*}
We have three conceivable cases as follows.\\
\begin{itemize}
\item If $h<p+q$, then we must have $u<q$ or $h-u-s<p$, otherwise $s<0$ which is impossible. Hence, $$g(U_{l_s}\cdots U_{l_1}e_{i_{j_r}\cdots i_{j_1}}t^{(u)})=0$$ or $$\dfrac{f(e_{i_1\cdots \widehat{i}_{j_1}\cdots \widehat{i}_{j_r}\cdots i_m}\wedge e_{l_1\cdots l_s}t^{(h-u-s)})}{x_{l_1}\cdots x_{l_s}}=0$$ and so is their product.
\item In case $h=p+q$, $g(U_{l_s}\cdots U_{l_1}\cdot e_{i_{j_r}\cdots i_{j_1}}t^{(u)})\cdot \dfrac{f(e_{i_1\cdots \widehat{i}_{j_1}\cdots \widehat{i}_{j_r}\cdots i_m}\wedge e_{l_1\cdots l_s}t^{(h-u-s)})}{x_{l_1}\cdots x_{l_s}} \neq 0$ only if $u=q$ and $h-u-s=p$. This implies that $s=0$. Then we have $$\{i_{j_1}, \ldots, i_{j_r}\}=J$$ and $$K\smallsetminus \{i_{j_1}, \ldots, i_{j_r}\}=I.$$ If $I \cap J = \emptyset$, there is only one such $K=I\cup J$ and one gets $$(f\smile  g)(e_Kt^{(h)})=\mathrm{sgn}(I,J)\mathbf{x}^{\alpha+\beta }.$$ If $I \cap J \neq \emptyset$, one has $K\smallsetminus J \subsetneq I$ for all $K$ which yields that $(f\smile  g)(e_Kt^{(h)})= 0$.
\item If $h>p+q$, then $s>0$. Indeed, by the argument as in above case, one has $u=q$ and $h-u-s=p$. Therefore, $s=h-p-q>0$. By the definition of the $U_{l_\bullet}$'s in Lemma \ref{lem2.chainmap} we can rewrite $U_{l_s}\cdots U_{l_1}\cdot e_{i_{j_r}\cdots i_{j_1}}t^{(u)}$ as a sum of some elements in the form below:
\begin{align*}
(x_1\cdots x_{z-1}\otimes x_{z+1}\cdots x_n)e_z\cdot a\otimes b \cdot e_Rt^{(u)}
\end{align*}
for some $z \in [n]$ and $a,b \in A$.\\
Applying the function $g$ on this sum, we get the result in $A$:
\begin{align*}
x_1\cdots x_{z-1}\cdot x_{z+1}\cdots x_n\cdot ab\cdot g(e_z\wedge e_Rt^{(u)}).
\end{align*}
Notice that as $g$ is a cocycle, by Lemma \ref{lem1.chain_map}, $x_z$ is a divisor of $g(e_z\wedge e_Rt^{(u)})$. Hence, $x_1\cdots x_{z-1}\cdot x_{z+1}\cdots x_n\cdot g(e_z\wedge e_Rt^{(u)})$ is a multiple of $x_1\cdots x_n$, which is zero in $A$.
\end{itemize}
Thus the result follows.
\qed
\end{Prop}
\section{The ring structure of $\mathrm{HH}^*(A)$}\label{s.5}
In this section, we are going to give a presentation for the algebra $\mathrm{HH}^*(A)$ by generators and relations.
\begin{Thm} \label{main.theorem}
Let $k[x_1,\ldots , x_n]$ be a polynomial ring over a field $k$ and let $A$ be the quotient ring $k[x_1,\ldots , x_n]/\langle x_1 \cdots x_n \rangle$. Then we have the isomorphism:
\begin{equation*}
HH^*(A)\cong k[X_1,\ldots , X_n, Y_1,\ldots , Y_n,Z]/\mathcal{I},
\end{equation*}
where $k[X_1,\ldots , X_n, Y_1,\ldots , Y_n,Z]$ is a graded commutative polynomial ring; $\mathrm{deg}\text{ } X_i=0$, $\mathrm{deg}\text{ } Y_i=1$ for all $i \in [n]$ and $\mathrm{deg}\text{ } Z=2$; and the ideal $\mathcal{I}$ is generated by the following relations:
\begin{itemize}
\item $a_1\cdots a_n \text{ where } a_i \in \{X_i,Y_i\};$
\item $ Y_i^2 \text{ for all } i \in [n]$;
\item $
\dfrac{X_1\cdots X_n}{{X}_{i_1}\cdots {X}_{i_m}}.\left(\sum\limits_{j=1}^m{(-1)^{j+1}Y_{i_1}\cdots \widehat{Y}_{i_j}\cdots Y_{i_m}}\right)Z
$, where $ m \in [n] \text{ and } 1\leq i_1 <i_2<\cdots <i_m\leq n.$
\end{itemize}
\proof
By the construction of Hochschild cohomology, $\mathrm{HH}^m(A)$ consists of the cosets of the cocycles of degree $m$. From the formula of multiplication, a cocycle of degree $m$ can be factorized into elements of degree 0, 1 and 2. Indeed, assume $E=\left(e_It^{(p)},\mathbf{x}^\alpha\right)$ is a cocycle where $I=\{i_1,\ldots ,i_s\}$, $\alpha=(\alpha_1,\ldots,\alpha_n)$. Let us write $x_1^{\alpha_1}\cdots x_n^{\alpha_n}$ as $x_{i_1}^{\alpha_{i_1}}\cdots x_{i_s}^{\alpha_{i_s}}\cdot x_{i_{s+1}}^{\alpha_{i_{s+1}}}\cdots x_{i_n}^{\alpha_{i_n}}$. Since $\partial\left(e_It^{(p)},\mathbf{x}^\alpha\right)=0$, it follows by Corollary \ref{cor1} $I \subseteq \mathrm{supp}(\mathbf{x}^\alpha)=\{i\mid \alpha_i > 0\}$, which means that $\alpha_{i_j} >0$ for all $j \in [s]$. Therefore, by the formula of the multiplication for two cochains in Proposition \ref{th.multiplication}, we obtain that:
\begin{align*}
E&=\left(e_{i_1}\cdots e_{i_s},x_{i_1}\cdots x_{i_s}\right)\cdot \left(t,1\right)^q\cdot \left(1,x_{i_1}^{\alpha_{i_1}-1}\cdots x_{i_s}^{\alpha_{i_s}-1}\cdot x_{i_{s+1}}^{\alpha_{i_{s+1}}}\cdots x_{i_n}^{\alpha_{i_n}}\right)\\
&=\left(e_{i_1},x_{i_1}\right)\cdots \left(e_{i_s},x_{i_s}\right) \left(t,1\right)^q \left(1,x_{i_1}\right)^{\alpha_{i_1}-1} \cdots \left(1,x_{i_s}\right)^{\alpha_{i_s}-1} \left(1,x_{i_{s+1}}\right)^{\alpha_{i_{s+1}}}\cdots \left(1,x_{i_n}\right)^{\alpha_{i_n}}.
\end{align*}
Briefly $E$ is factorized into the following elements:
$\left(t,1\right)$, $q$ times; $\left(e_i,x_i\right)$ where $i \in I$; and $\left(1,x_j\right)$, $\alpha_j-\beta_j$ times, where $j \in \mathrm{supp}(\mathbf{x}^{\alpha})$, $\beta_j=1$ if $j \in I$ and $\beta=0$ if $j \notin I$.\\
For each $i\in [n]$, set $X_i$ to be the coset of the element $\left(1,x_i\right)$, $Y_i$ to be the coset of the element  $\left(e_i,x_i\right)$ and $Z$ to be the coset of the element  $\left(t,1\right)$. Then we have $\mathrm{HH}^*(A)$ is generated by  $X_i$, $Y_i$ and $Z$. As $x_1\cdots x_n=0$, we obtain the relations $
a_1\cdots a_n \text{ where } a_i \in \{X_i,Y_i\}
$. The relations $Y_i^2$ come from the fact that $e_i^2=0$. \\
\begin{Rem} \label{rem1}
For any element $\left(e_It^{(p)},\mathbf{x}^\alpha\right)$, by Lemma \ref{lemma1} we obtain that the image $\partial \left(e_It^{(p)},\mathbf{x}^\alpha\right)$ is a multiple of $\partial \left(e_I,1\right)$. Hence, another relation is $\partial \left(e_I,1\right)$. Suppose that $I=\lbrace i_1,i_2,\ldots ,i_m\rbrace$ where $m \in [n]$ and $1 \leq i_1<i_2<\cdots <i_m\leq n$. It follows that
\begin{equation*}
\partial \left(e_I,1\right)= \sum_{j=1}^{m}{(-1)^{j+1}\left(e_{i\cdots \widehat{i}_j\cdots i_m}t,\frac{x_1x_2\cdots x_n}{x_{i_j}}\right)}
\end{equation*}
By relabeling the indices in $x_1x_2\cdots x_n$ as $x_{i_1}x_{i_2}\cdots x_{i_m}\cdot x_{i_{m+1}}\cdots x_{i_n}$, we rewrite $\partial \left(e_I,1\right)$ as a combination of generators as follows:
\begin{align*}
\partial \left(e_I,1\right)=&\sum_{j=1}^{m}{(-1)^{j+1}\left(e_{i\cdots \widehat{i}_j\cdots i_m}t,x_{i_1}\cdots \widehat{x}_{i_j}\cdots x_{i_m}\cdot x_{i_{m+1}}\cdots x_{i_n}\right)}
\\
=&\left(\sum_{j=1}^{m}(-1)^{j+1}(e_{i_1},x_{i_1})\cdots \widehat{(e_{i_j},x_{i_j})}\cdots (e_{i_m},x_{i_m})\right)(t,1)(1,x_{i_{m+1}})\cdots (1,x_{i_n}).
\end{align*}
\end{Rem}
We have shown that $X_1,\ldots , X_n, Y_1,\ldots , Y_n,Z$ generate $\mathrm{HH}^*(A)$ and that they satisfy the relations in $\mathcal{I}$. Now we prove that there exists the isomorphism as in the assertion.\\

Let $S:= k[X_1,\ldots , X_n, Y_1,\ldots , Y_n,Z]$ be the graded commutative polynomial ring over the field $k$ and $\mathcal{J}$ be the ideal of $S$ generated by the elements $a_1\cdots a_n \text{ where } a_i \in \{X_i,Y_i\}$, and $Y_i^2 \text{ for all } i \in [n]$. As $\mathcal{J}$ is a monomial ideal, the residue classes of the monomials not belonging to $\mathcal{J}$ form a $k$-basis of the quotient ring $S/\mathcal{J}$. The monomial $X_1^{\alpha_1}\cdots  X_n^{\alpha_n} \cdot Y_1^{\beta_1}\cdots  Y_n^{\beta_n}\cdot Z^{q} \in S$ is not in $\mathcal{J}$ if and only if $\beta_i \leq 1$ for all $i$ and $\{i\mid \alpha_i >0 \text{ or } \beta_i>0\}\subsetneq [n]$. We can identify a non-zero residue class in $S/\mathcal{J}$ by the $k$-basis element which represents it. Let us construct the map $\psi$ from $S/\mathcal{J}$ to $\mathrm{Ker}(\partial)$ by sending the $k$-basis element $X_1^{\alpha_1}\cdots  X_n^{\alpha_n} \cdot Y_1^{\beta_1}\cdots  Y_n^{\beta_n}\cdot Z^{q}$ in $S$ to $(e_1^{\beta_1}\cdots e_n^{\beta_n}t^{(q)},x_1^{\alpha_1+\beta_1}\cdots x_n^{\alpha_n+\beta_n})$ in $\mathrm{Ker}(\partial)$. We can check that $\psi$ is an isomorphism between these algebras. \\
By Remark \ref{rem1}, the image $\mathrm{Im}(\partial)$ is generated by the relations $\partial \left(e_I,1\right)$ and $$\psi \left(
\dfrac{X_1\cdots X_n}{{X}_{i_1}\cdots {X}_{i_m}}.\left(\sum\limits_{j=1}^m{(-1)^{j+1}Y_{i_1}\cdots \widehat{Y}_{i_j}\cdots Y_{i_m}}\right)Z\right)=\partial \left(e_I,1\right),
$$
where  $I=\lbrace i_1,i_2,\ldots ,i_m\rbrace$. Therefore, $\psi^{-1}(\mathrm{Im}(\partial))=\mathcal{I}$. Hence, $S/\mathcal{I}\cong \dfrac{\mathrm{Ker}(\partial)}{\mathrm{Im}(\partial)}=\mathrm{HH}^*(A)$.
\qed
\end{Thm}
\begin{Ex}
Applying the above result to the case $n=2$, we have that:
$$\mathrm{HH}^*(k[x,y]/\langle xy\rangle) \cong k[x_1,x_2,y_1,y_2,z]/\mathcal{I}$$
where $\mathrm{deg}\text{ }x_i=0$, $\mathrm{deg}\text{ }y_i=1$ for $i=1,2$, $\mathrm{deg}\text{ }z=2$ and the ideal $\mathcal{I}$ is generated by $x_1x_2$, $x_1y_2$, $y_1x_2$, $y_1y_2$, $y_1^2$, $y_2^2$, $x_1z$, $x_2z$, $(y_1+y_2)z.$
\end{Ex}
\section{The Hilbert series of $\mathrm{HH}^*(A)$}\label{s.6}
In this final section, we apply the results to compute the Hilbert series of $\mathrm{HH}^*(A)$. First we introduce a grading on $\mathrm{HH}^*(A)$, which is combined from two different gradings. The first is the $\mathbb{N}$-grading based on the degree of homology. In detail, if $(e_It^{(q)},\mathbf{x}^\alpha)$ is an element in $\mathrm{Hom}_k(\overline{F}_m,A)$ (which means $|I|+2q=m$), we let $\mathrm{hdeg}(e_It^{(q)},\mathbf{x}^\alpha)=m$. The other one is the $\mathbb{Z}^n$-grading based on the lattice point representation for a variable. Let us set $\mathrm{rdeg}(t,1)=-(1,1,\ldots ,1)$ (an $n$-vector with all 1s), $\mathrm{rdeg}(e_i,1)=-\mathbf{e}_i$ and $\mathrm{rdeg}(1,x_i)=\mathbf{e}_i$ for all $i \in \{1,2,\ldots ,n\}$, where $\mathbf{e}_i$ is the $i$th standard basis vector in $\mathbb{N}^n$. The differential $\partial$ is a 1-homogeneous morphism with respect to the first grading and a 0-homogeneous morphism with respect to the second grading. We call the grading given by combining these gradings the \textit{multidegree} of a standard element, $\mathrm{mdeg}(e_It^{(q)},\mathbf{x}^\alpha):=(\mathrm{hdeg}(e_It^{(q)},\mathbf{x}^\alpha),\mathrm{rdeg}(e_It^{(q)},\mathbf{x}^\alpha))$ in $\mathbb{N}\times \mathbb{Z}^n$. Equivalently, if the element $(e_It^{(q)},\mathbf{x}^\alpha)$ whose component $e_I$ is identified by the vector $(\epsilon_1,\epsilon_2,\ldots ,\epsilon_n)$ where for any $i$ in $\{1,2,\ldots,n\}$, $\epsilon_i$ is 1 if $i \in I$ and 0 otherwise and $\mathbf{x}^\alpha$ is identified with $\alpha=(\alpha_1,\alpha_2,\ldots ,\alpha_n)$, then we have
\begin{equation*}
\mathrm{mdeg}(e_It^{(q)},\mathbf{x}^\alpha)=(|I|+2q,\alpha_1-\epsilon_1-q,\ldots,\alpha_n-\epsilon_n-q)
\end{equation*}
Thus the element $(e_It^{(q)},\mathbf{x}^\alpha)$ contributes the term
\begin{equation*}
a_0^{|I|+2q}a_1^{\alpha_1-\epsilon_1-q}a_2^{\alpha_2-\epsilon_2-q}\cdots a_n^{\alpha_n-\epsilon_n-q}
\end{equation*}
(or briefly, as $\mathbf{a}^{\chi}$ where $\chi=\mathrm{mdeg}(e_It^{(q)},\mathbf{x}^\alpha)$) to the Hilbert series.\\
Let $\mathbf{H}_\chi$ be the $k$-module generated by the elements whose multidegree is $\chi \in \mathbb{N}\times\mathbb{Z}^n$. The Hilbert series of $\mathrm{HH}^*(A)=\oplus_{\chi \in \mathbb{N}\times\mathbb{Z}^n}\mathbf{H}_{\chi}$ as an  $\mathbb{N}\times \mathbb{Z}^n$-graded vector space via the grading above is the formal power series:
\begin{equation*}
\mathcal{H}(\mathrm{HH}^*(A);\mathbf{a})=\sum_{\chi \in \mathbb{N}\times\mathbb{Z}^n}\mathrm{dim}_k(\mathbf{H}_{\chi})\mathbf{a}^{\chi}.
\end{equation*}
\begin{Thm} \label{hilbert.theorem}
The Hochschild cohomology ring $\mathrm{HH}^*(A)$ has the Hilbert series:
\begin{equation*}
\mathcal{H}(\mathrm{HH}^*(A);\mathbf{a})=\dfrac{(a_0+1)^{n+1}a_1\cdots a_n-(a_0+a_1)\cdots (a_0+a_n)\cdot\left(a_0+a_1\cdots a_n\right)}{(a_1\cdots a_n-a_0^2)\cdot(1-a_1)\cdots (1-a_n)}.
\end{equation*}
\proof
Let us denote $\mathrm{H}_1$ and $\mathrm{H}_2$ the Hilbert series of the cocycles and the coboundaries respectively. As $|I|=\epsilon_1+\ldots+\epsilon_n$, it follows by a simple computation that:
\begin{equation*}
a_0^{|I|+2q}\cdot a_1^{\alpha_1-\epsilon_1-q}\cdots a_n^{\alpha_n-\epsilon_n-q}=(a_0^2(a_1\cdots a_n)^{-1})^q\cdot \prod\limits_{i = 1}^n {{{(a_0a_i^{ - 1})}^{{\epsilon _i}}}a_i^{{\alpha _i}}}.
\end{equation*}
From Corollary \ref{cor1.2}, we recall that the standard element $(e_It^{(q)},\mathbf{x}^\alpha)$ is a cocycle if $I \subseteq \mathrm{supp}(\mathbf{x}^\alpha)$, i.e., for any $i \in [n]$ we have $\alpha_i>0$ if $\epsilon_i=1$ and $\alpha_i\geq 0$ if $\epsilon_i=0$. Consequently, we need to eliminate the cases that all $\alpha_i> 0$.
Thus we get the first series which counts all cocycles:
\begin{equation*}
\mathrm{H}_1=\dfrac{1}{1-a_0^2(a_1 \cdots a_n)^{-1}}\cdot \left[\prod\limits_{i = 1}^n {\dfrac{a_0+1}{1-a_i}}- \prod\limits_{i=1}^n{\dfrac{a_0+a_i}{1-a_i}}\right].
\end{equation*}
To obtain the series of the coboundaries $\mathrm{H}_2$, we consider the element of the form $(e_It^{(q)},\mathbf{x}^\alpha)$ in the same multidegree as above. Then the multidegree of the image $\partial(e_It^{(q)},\mathbf{x}^\alpha$) is
\begin{equation*}
(|I|+2q+1,\alpha_1-\epsilon_1-q,\ldots,\alpha_n-\epsilon_n-q).
\end{equation*}
All cases, $|I\setminus \mathrm{supp}(\mathbf{x}^\alpha)|=m$ for $m$ from 1 to $n$, are counted, except for $m=0$ (which means $I \subseteq \mathrm{supp}(\mathbf{x}^\alpha)$ and hence, $\partial(e_It^{(q)},\mathbf{x}^\alpha)=0$). Then we get the series of the coboundaries:
\begin{equation*}
\mathrm{H}_2=\dfrac{a_0}{1-a_0^2(a_1 \cdots a_n)^{-1}}\cdot \left[\prod\limits_{i = 1}^n {\dfrac{a_0a_i^{-1}+1}{1-a_i}}- \prod\limits_{i=1}^n{\dfrac{a_0+1}{1-a_i}}\right].
\end{equation*}
Now we are able to get the Hilbert series of $\mathrm{HH}^*(A)$, which is $\mathrm{H}_1-\mathrm{H}_2$.
\qed
\end{Thm}

\end{document}